\documentclass[12pt]{article}

\usepackage{amsmath,amssymb,amsthm}
\usepackage{cite}

\newtheorem{theorem}{Theorem}[section]
\newtheorem{corollary}[theorem]{Corollary}
\newtheorem{conjecture}[theorem]{Conjecture}
\newtheorem{lemma}[theorem]{Lemma}

\newcommand{\address}[1]{#1}

\DeclareMathOperator{\cp}{\,\square\,}

\begin{document}

\title{On decycling and forest numbers of Cartesian products of trees}
\author{
Ali Ghalavand $^{a}$ \\
\texttt{\footnotesize alighalavand@nankai.edu.cn}  \and Sandi Klav\v{z}ar $^{b,c,d}$ \\ \texttt{\footnotesize sandi.klavzar@fmf.uni-lj.si}
\and Ning Yang $^{a,}$\footnote{Corresponding author} \\ \texttt{\footnotesize yn@mail.nankai.edu.cn}}

\maketitle	

\address{
\noindent
    $^a$ Center for Combinatorics and LPMC, Nankai University, Tianjin 300071, China\\
    $^b$ Faculty of Mathematics and Physics, University of Ljubljana, Slovenia\\
	$^c$ Institute of Mathematics, Physics and Mechanics, Ljubljana, Slovenia\\
	$^d$ Faculty of Natural Sciences and Mathematics, University of Maribor, Slovenia
}

\begin{abstract}
The decycling number $\nabla(G)$ of a graph $G$ is the minimum number of vertices that must be removed to eliminate all cycles in $G$. The forest number $f(G)$ is the maximum number of vertices that induce a forest in $G$. So $\nabla(G) + f(G) = |V(G)|$. For the Cartesian product $T \,\square\, T'$ of trees $T$ and $T'$ it is proved that $\nabla(S_n \,\square\, S_{n'}) \leq \nabla(T \,\square\, T')$, thus resolving the conjecture of Wang and Wu asserting that $f(T \,\square\, T') \leq f(S_n \,\square\, S_{n'})$. It is shown that $\nabla(T \,\square\, T') \ge |V(T)| - 1$ and the equality cases characterized. For prisms over trees, it is proved that $\nabla(T\,\square\, K_2) = \alpha'(T)$, and  for arbitrary graphs $G_1$ and $G_2$, it is proved that $\nabla(G_1 \,\square\, G_2) \geq \alpha'(G_1) \alpha'(G_2)$, where $\alpha'$ is the matching number.
\end{abstract}


\medskip\noindent
{\bf AMS Subj.\ Class.\ (2020)}: 05C35, 05C69, 05C76

%%%%%%%%%%%%%%%%%%%%%%%%%%%%%%%%%%
\section{Introduction}
%%%%%%%%%%%%%%%%%%%%%%%%%%%%%%%%%%

A {\em decycling set} \(S\) of a graph \(G\) is a subset of its vertices such that \(V(G) - S\) is a forest. In other words, each cycle of the graph contains at least one vertex of a decycling set. From this reason, decycling sets are also known in the literature as feedback vertex sets, cf.~\cite{dabrowski-2024, dekker-2024}. The {\em decycling number} \(\nabla(G)\) of $G$ is the minimum cardinality of a decycling set. The {\em forest number} \(f(G)\) of $G$ is the maximum number of vertices that induce a forest in \(G\). The following equation that goes back to Erd\"{o}s  et al.~\cite{Erdos1} is evident:
\begin{equation}
\label{eq:key-connection}
\nabla(G) + f(G) = |V(G)|\,.
\end{equation}

The issue of eliminating all cycles in a graph by deleting a set of vertices was first introduced in 1974 in combinatorial circuit design~\cite{j1}. Soon after it has been established that determining the decycling number is NP-complete for general graphs~\cite{kar1}.
Decycling sets have found applications across various fields, including deadlock prevention in operating systems~\cite{si1, w1}, constraint satisfaction problem and Bayesian inference in artificial intelligence~\cite{ba1}, monopolies in synchronous distributed systems~\cite{pe1, pe2}, the converter placement problem in optical networks~\cite{kl1}, and VLSI chip design by~\cite{fe1}.

At this point, we should also mention the related concept of cycle isolation, where we also want to destroy all cycles of a given graph, the difference here is that together with a given set of vertices we remove its entire closed neighbourhood. For more information on cycle isolation see~\cite{borg-2020, liang-2026, zgang-2024} and references therein.

In 1987, Akiyama and Watanabe \cite{akiyama1}, along with Albertson and Haas in 1998 \cite{albertson1}, independently conjectured that if $G$ is a planar bipartite graph, then $f(G) \ge \frac{5}{8}|V(G)|$. Motivated by this conjecture, Alon~\cite{alon1} proved that every bipartite graph $G$ with an average degree \( d \geq 1 \) has $f(G) \ge \left(\frac{1}{2} + e^{-bd^2}\right)|V(G)|$ for some positive constant \( b \). On the other hand, there are bipartite graphs $G$ with an average degree \( d \geq 1 \) such that $f(G) < \left(\frac{1}{2} + \frac{1}{e^{b'\sqrt{d}}}\right)|V(G)|$, where $b'$ is an absolute constant. In 2014, Conlon et al.~\cite{conlon1} improved Alon's lower bound to \( \left(\frac{1}{2} + d^{-db}\right)|V(G)| \) for \( d \geq 2 \).
In 2017, Wang et al. \cite{wang2} demonstrated that every simple bipartite planar graph \( G \) satisfies the inequality \( f(G) \ge \lceil \frac{4|V(G)| + 3}{7} \rceil \). More recently, one of the authors, along with other \cite{Ghalavand1}, proved that every balanced bipartite graph \( G \) with a minimum degree of at least \( \frac{|V(G)|}{4} + 1 \) has \( f(G) = \frac{|V(G)|}{2} + 1 \). This result was initially proposed as a conjecture in \cite{wang1}.
 Upper and lower bounds for the forest number of a graph in relation to its order, size, and maximum degree are presented in~\cite{beineke1,shi1}. The forest number has been further studied in the context of planar graphs \cite{dross1, dross2, kelly1, le1, ma-2024, Petrusevski1}, regular graphs \cite{punnim1,ren1}, (sub)cubic graphs \cite{kelly2, nedela}, line graphs~\cite{ma-2025}, and graphs  embedded into surfaces~\cite{ma-2025b}. 

The notations \( P_n \),  \( S_n \), and \(C_n\) represent the path, the star, and the cycle graph on \( n \) vertices, respectively.

The study of decycling and forest numbers of Cartesian products of graphs has a long history. Beineke and Vandell~\cite{beineke1} investigated them on grids and hypercubes. Luccio~\cite{n-Luccio1} demonstrated that if \( n, n' \geq 2 \), then \( \nabla(P_n \cp P_{n'}) \geq \left\lceil \frac{(n-1)(n'-1)+1}{3} \right\rceil \). Madelaine and Stewart~\cite{Madelaine-Stewart} established that \( \nabla(P_n \cp P_{n'}) \leq \left\lceil \frac{(n-1)(n'-1)+1}{3} \right\rceil + 2 \). Lien et al.~\cite{Lien-Fu-Shih} followed by proving \( \nabla(P_n \cp P_{n'}) \leq \left\lceil \frac{(n-1)(n'-1)+1}{3}  \right\rceil+ 1 \). Pike and Zou~\cite{Pike-Zou} derived \( \nabla(C_n \cp C_{n'}) = \left\lceil \frac{3n'}{2} \right\rceil \) when \( n = 4 \), and \( \nabla(C_n \cp C_{n'}) = \left\lceil \frac{nn'+2}{3} \right\rceil \) otherwise. Recently, and as our main motivation, Wang and Wu~\cite{wang1} investigated the forest number of Cartesian, direct, and lexicographic products of graphs. Among other things, in \cite[Conjecture 5.1]{wang1} they put forward the following conjecture, where $S_n$ denotes the star $K_{1,n-1}$.

\begin{conjecture}
\label{conj1}
If $T$ and $T'$ are trees of respective orders $n$ and $n'$, then
\[f(T \cp T') \leq f(S_{n} \cp S_{n'})\,.\]
\end{conjecture}

In Section~\ref{sec:conjecture-upper}, we prove this conjecture. In the subsequent section, we characterize the trees that attain the equality in the conjecture. If the orders of $T$ and $T'$ are the same, then the equality is attained only if at least  one of the factors is a star, while if the orders of $T$ and $T'$ are different, then the equality is attained if and only if the larger factor is a star (and the smaller factor is arbitrary).  In Section~\ref{sec:more}, we consider the missing case when a star is the smaller factor and establish a sharp lower bound for the decycling number of Cartesian products of graphs based on the matching number of their factors.

To conclude the introduction, we give definition needed in this paper.

 The subgraph of a graph $G$ induced by \(X\subseteq V(G)\) is denoted by \(G[X]\).  A {\em vertex cover} of \(G\) is a subset of its vertices that includes at least one endpoint of every edge. The {\em covering number} \(\beta(G)\) of \(G\) is the minimum size of a vertex cover of $G$. A {\em matching} of \(G\) is a subset of edges such that no two edges from it share a common vertex. The {\em matching number}  \(\alpha'(G)\) of \(G\) is the maximum size of a matching set for the graph.

Let \(G\) and \(H\) be two graphs. Their {\em Cartesian product} \(G \cp H\) is a graph with vertex set \(V(G) \times V(H)\), and two vertices \((g, h)\) and \((g', h')\) are adjacent if either \(g = g'\) and \(hh' \in E(H)\), or \(gg' \in E(G)\) and \(h = h'\). The subgraph of $G\cp H$ induced by the set of vertices $\{(g,h):\ g\in V(G)\}$ is called a {$G$-layer} and denoted by $G^h$. Similarly, the subgraph of $G\cp H$ induced by $\{(g,h):\ h\in V(H)\}$ is an {\em $H$-layer} and denoted by $^gH$. Note that $G^h\cong G$ for every $g\in V(G)$, and that $^gH \cong H$ for every $h\in V(H)$. Recall that the Cartesian product operation is commutative, that is, $G\cp H \cong H\cp G$, where $\cong$ denotes graph isomorphism.

%%%%%%%%%%%%%%%%%%%%%%%%%%%%%%%%%
\section{Proof of Conjecture~\ref{conj1}}
\label{sec:conjecture-upper}
%%%%%%%%%%%%%%%%%%%%%%%%%%%%%%%%%

Conjecture~\ref{conj1} will follow from the following result.

\begin{theorem}\label{thm:main}
If $T$ and $T'$ are trees of respective orders $n$ and $n'$, where $n'\geq n\geq 2$, then
\[\nabla(T \cp T')\geq n-1 = \nabla(T \cp S_{n'}).\]
\end{theorem}

\proof
Let $T$ and $T'$ be trees $V(T)=\{v_1,\ldots,v_n\}$ and $V(T')=\{u_1,\ldots,u_{n'}\}$, where $n'\geq n\geq 2$. Suppose on the contrary that $\nabla(T \cp T')\leq n-2$ and let $S$ be a decycling set in $T \cp T'$ with $|S|=\nabla(T \cp T') \le n-2$. By the pigeonhole principle there exist $i, j\in [n]$, $i\ne j$, such that $S\cap V(^{v_i}T') = \emptyset$ and $S\cap V(^{v_j}T') = \emptyset$, and there exist $k, \ell\in [n']$, $k\ne \ell$, such that $S\cap V(T^{u_k}) = \emptyset$ and $S\cap V(T^{u_\ell}) = \emptyset$. Consider now the subgraph $X$ of $T \cp T'$ induced by
$$V(^{v_i}T') \cup V(^{v_j}T') \cup V(T^{u_k}) \cup V(T^{u_\ell})\,.$$
$X$ is a connected graph of order $2n + 2n' - 4$ and of size at least $2n + 2n' - 4$. It follows that $X$ contains a cycle which is not possible because $S\cap V(X) = \emptyset$. This contradiction proves that $\nabla(T \cp T')\geq n-1$.

By the above, $\nabla(T \cp S_{n'}) \ge n-1$. Let $V(S_{n'})=\{u_1,\ldots,u_{n'}\}$, and let $\{u_1,\ldots,u_{n'-1}\}$ be an independent set of $S_{n'}$. In this case, the  subgraph of $T \cp S_{n'}$ induced by $V(T)\times (V(S_{n'})-\{u_{n'}\})$ is isomorphic to the disjoint union of $(n'-1)$ copies of $T$. Consequently, the subgraph induced by
$$\left(V(T)\times (V(S_{n'})-\{u_{n'}\})\right)\cup\{(v_{1},u_{n'})\}$$
is a tree. It follows that $\{(v_{i},u_{n'}):\ 2 \leq i \leq n\}$ is a decycling set for $T \cp S_{n'}$, and thus $ \nabla(T \cp S_{n'}) \leq n-1$.
\qed

\medskip
Theorem~\ref{thm:main} in particular implies that $\nabla(T \cp T')\geq n-1 = \nabla(S_n \cp S_{n'})$. Conjecture~\ref{conj1} now follows by applying~\eqref{eq:key-connection}. More precisely, since $\nabla(T \cp T')\geq \nabla(S_n \cp S_{n'})$, by~\eqref{eq:key-connection} we get $|V(T \cp T')|-\nabla(T \cp T')\leq |V(T \cp T')| - \nabla(S_n \cp S_{n'})$. Thus we can conclude that
\begin{align*}
f(T \cp T') & = |V(T \cp T')| - \nabla(T \cp T') \\
& \leq |V(T \cp T')| - \nabla(S_n \cp S_{n'}) \\
& = f(S_n \cp S_{n'})\,,
\end{align*}
as desired.

\medskip
Using the same line of reasoning as above, Theorem~\ref{thm:main} also yields:

\begin{corollary}
If $T$ is a nontrivial tree, and $n'\ge |V(T)|$ an integer, then
\[f(T\cp S_{n'}) = |V(T)|n'-|V(T)| + 1\,.\]
\end{corollary}

%%%%%%%%%%%%%%%%%%%%%%%%%%%%%%%%%
\section{The equality case in Theorem~\ref{thm:main}}
\label{sec:equality}
%%%%%%%%%%%%%%%%%%%%%%%%%%%%%%%%%

In this section, we characterize the graphs that attain the equality in Theorem~\ref{thm:main} and consequently in Conjecture~\ref{conj1}. For this sake, the following two results are crucial.

\begin{lemma}\label{lem:at least n first}
If $T$ and $T'$ are trees of respective orders $n$ and $n'$, where $n'> n\geq 2$, and $T'\not\cong S_{n'}$, then
\[\nabla(T \cp T')\geq n\,.\]
\end{lemma}

\proof
Let $V(T)=\{v_1,\ldots,v_n\}$ and $V(T')=\{u_1,\ldots,u_{n'}\}$. Suppose on the contrary that $\nabla(T \cp T') \leq n-1$ and let $S$ be a decycling set in $T \cp T'$ with $|S|=\nabla(T \cp T')\le n-1$. We distinguish the following two cases.

\medskip\noindent
{\bf Case 1}: $\exists\ i, j\in [n]$, $i\ne j$, such that $S\cap V(^{v_i}T') = \emptyset$ and $S\cap V(^{v_j}T') = \emptyset$.\\
Since $n'>n$, there also exist $k, \ell\in [n']$, $k\ne \ell$, such that $S\cap V(T^{u_k}) = \emptyset$ and $S\cap V(T^{u_\ell}) = \emptyset$. As in the proof of Theorem~\ref{thm:main} we now consider the subgraph of $T \cp T'$ induced by $V(^{v_i}T') \cup V(^{v_j}T') \cup V(T^{u_k}) \cup V(T^{u_\ell})$ to arrive at a contradiction.

\medskip\noindent
{\bf Case 2}: $S\cap V(^{v_i}T') = \emptyset$, and $|S\cap V(^{v_j}T')| = 1$ for every $j\in [n]$, $j\ne i$. \\
Since $T'$ is not a star, it contains an induced path on four vertices, say $u_{j_1}, u_{j_2}, u_{j_3}, u_{j_4}$. Let further $v_\ell$ be a neighbor of $v_i$ in $T$. Let $Y$ be the subgraph of
$T\cp T'$ induced by the vertices
$$ \{v_i, v_\ell\} \times \{u_{j_1}, u_{j_2}, u_{j_3}, u_{j_4}\}\,.$$
The subgraph $Y$ is isomorphic to $P_2 \cp P_4$ and since $|S\cap V(Y)| \le 1$ we can conclude that $S$ is not a decycling set.
\qed

\begin{lemma}\label{lem:at least n second}
If $T$ and $T'$ are trees, both of order $n \geq 4$, and none of $T$ and $T'$ is a star, then
\[\nabla(T \cp T') \geq n.\]
\end{lemma}

\proof
Let $V(T)=\{v_1,\ldots,v_n\}$, $V(T')=\{u_1,\ldots,u_n\}$, and suppose for the sake of contradiction that $\nabla(T \cp T') \leq n-1$. Let $S$ be decycling set in $T \cp T'$ with $|S|=\nabla(T \cp T') \le n-1$. If there exist $i, j\in [n]$, $i\ne j$, such that $S\cap V(^{v_i}T') = \emptyset$ and $S\cap V(^{v_j}T') = \emptyset$, and there exist $k, \ell\in [n]$, $k\ne \ell$, such that $S\cap V(T^{u_k}) = \emptyset$ and $S\cap V(T^{u_\ell}) = \emptyset$, then we get a contradiction as in the proof of Theorem~\ref{thm:main}.

By the commutativity of the Cartesian product we may without loss of generality assume that in the remaining case there exists an integer $i\in [n]$ such that $S\cap V(^{v_i}T') = \emptyset$, and $|S\cap V(^{v_j}T')| = 1$ for every $j\in [n]$, $j\ne i$. Then we get a contradiction along the same lines as in in Case 2 of the proof of Lemma~\ref{lem:at least n first}.
\qed

\medskip
The equality case in Theorem~\ref{thm:main} now reads as follows.

\begin{theorem}\label{thm:equality}
Let $T$ and $T'$ be trees of respective orders $n$ and $n'$, where $2\le n \leq n'$. Then the following holds.
\begin{enumerate}
\item[(i)] If $n=n'$, then $\nabla(T \cp T')=n-1$ if and only if at least one of $T\cong S_n$ and $T'\cong S_{n'}$ holds.
\item[(ii)] If $n < n'$, then $\nabla(T \cp T')=n-1$ if and only if $T'\cong S_{n'}$ $($and $T$ is arbitrary$)$.
\end{enumerate}
\end{theorem}

\proof
Combine Theorem~\ref{thm:main}, Lemmas~\ref{lem:at least n first} and~\ref{lem:at least n second}, and the facts $\nabla(P_2 \cp P_2)=1$ and $\nabla(P_3 \cp P_3)=2$.
\qed

\medskip
Using Theorem~\ref{thm:main}, Lemmas~\ref{lem:at least n first} and~\ref{lem:at least n second}, and~\eqref{eq:key-connection}, we also deduce the following.

\begin{corollary}\label{cor:equality-fof-f}
If $T$ and $T'$ are trees of respective orders $n$ and $n'$, where $2\le n \leq n'$, then $f(T \cp T') \leq nn'-n+1$. Moreover, the following holds.
\begin{enumerate}
\item[(i)] If $n=n'$, then $f(T \cp T')=nn'-n+1$ if and only if at least one of $T\cong S_n$ and $T'\cong S_{n'}$ holds.
\item[(ii)] If $n < n'$, then $f(T \cp T')=nn'-n+1$ if and only if $T'\cong S_{n'}$ $($and $T$ is arbitrary$)$.
\end{enumerate}
\end{corollary}

%%%%%%%%%%%%%%%%%%%%%%%%%%%%%%%%%
\section{Additional results}
\label{sec:more}
%%%%%%%%%%%%%%%%%%%%%%%%%%%%%%%%%

We have already proven that if \( T \) is a tree of order \( n\ge 2 \), and if $n'\ge n$, then \( \nabla(T\cp S_{n'}) = n - 1\). In the case when the star has order smaller than $T$, the following holds.

\begin{theorem}\label{thm:small-star}
If $2\le n' < n$, and $T\not\cong S_n$ is a tree of order $n$, then
\[n' \leq \nabla(T\cp S_{n'}) \leq n-1\,.\]
\end{theorem}

\proof
The lower bound follows by Lemma~\ref{lem:at least n first}. To establish the upper bound, let $V(T)=\{v_1,\ldots,v_n\}$ and $V(S_{n'}) = \{u_1,\ldots,u_{n'}\}$, where $u_{n'}$ is the central vertex of $S_{n'}$. Just as in the last part of the proof of Theorem~\ref{thm:main}, we deduce that the set $\{(v_i, u_{n'}):\ 2 \leq i \leq n\}$ is a decycling set in \(T \cp S_{n'}\). Consequently, \(\nabla(T\cp S_{n'}) \leq n-1\).
\qed

\begin{corollary}
If $T\not\cong S_{n+1}$ is a nontrivial tree of order $n+1$, then
  \[\nabla(S_{n} \cp T)=n\quad \text{and}\quad f(S_{n} \cp T)=n^2\,.\]
\end{corollary}

In the following theorem, we establish a relationship between the decycling number of Cartesian products and the matching number of their factors. For its proof recall that the famous K\"{o}nig-Egerv\'ary theorem~\cite{k-m, e-m} asserts that if $G$ is a bipartite graph, then  $\beta(G)=\alpha'(G)$.

\begin{theorem}\label{thm-general-matching}
If $G_1$ and $G_2$ are graphs, then
$$\nabla(G_1 \cp G_2) \geq \alpha'(G_1) \alpha'(G_2)\,.$$
Moreover, if $T$ is a tree, then $\nabla(T \cp P_2) = \alpha'(T)$.
\end{theorem}

\proof
To prove the first assertion, consider maximum matchings $M_1$ and $M_2$ of $G_1$ and $G_2$, respectively. Let $V_i$, $i\in [2]$, be the set of endvertices of the edges from $M_i$. Then the subgraph of $G_1 \cp G_2$ induced by the vertices $V_1\times V_2$ is isomorphic to the disjoint union of $\alpha'(G_1)\alpha'(G_2)$ copies of $C_4$. It follows that any decycling set $S$ in $G_1 \cp G_2$ must contain at least $\alpha'(G_1)\alpha'(G_2)$ elements, that is, $\nabla(G_1 \cp G_2)\geq\alpha'(G_1)\alpha'(G_2)$.

Let $T$ be a tree. Then by the already proved inequality we have $\nabla(T \cp P_2) \geq \alpha'(T)$. To prove that $\nabla(T \cp P_2) \leq \alpha'(T)$, consider a smallest vertex cover $W = \{w_1,\ldots,w_{\beta(T)}\}$ of $T$. By the definition of a vertex cover,  $T[V(T)-W]$ is an edgeless graph, and by K\"{o}nig-Egerv\'ary theorem, $|W|=\alpha'(T)$. Let $V(P_2)=\{a,b\}$ and set $\Omega=\{(a,w_i):\ i\in [\beta(T)]\}$. Since $T[V(T)-W]$ is an edgeless graph, the subgraph of $T \cp P_2$ induced by $V(T \cp P_2)\setminus \Omega$ is a tree. Therefore, $\nabla(T \cp P_1) \leq |W| = \alpha'(T)$.
\qed

%%%%%%%%%%%%%%%%%%%%%%%%%%%%%%%%%%%
\section{Concluding remarks}
\label{sec:conclude}
%%%%%%%%%%%%%%%%%%%%%%%%%%%%%%%%%%%

The main contribution of this paper is to demonstrate Conjecture~\ref{conj1} asserting that if $T$ and $T'$ are trees of respective orders $n$ and $n'$, then $f(T \cp T') \leq f(S_{n} \cp S_{n'})$, and to characterize graphs that attain the inequality. In~\cite[Conjecture 5.1]{wang1}, it was also conjectured that
$$f(P_{n} \cp P_{n'}) \leq f(T \cp T')$$
holds. This remains open. In partial support of the fact that this conjecture is difficult to verify, we note that  $\nabla(P_{n} \cp P_{n'})$ (equivalently $f(P_{n} \cp P_{n'})$) is not known yet in all the cases. The decycling number of the Cartesian product of paths has been investigated in~\cite{bau-2001}, for the summary of this research see Section 3 in the survey~\cite{bau-2002}. Later, progress on $\nabla(P_{n} \cp P_{n'})$ has been made in~\cite{Lien-Fu-Shih, n-Luccio1, Madelaine-Stewart}, however, to the best of our knowledge, the general solution is not yet known.

%%%%%%%%%%%%%%%%%%%%%%%%%%%%%%%%%%%%%%%%%%%%%%%%%%%%%%%
\section*{Acknowledgments}
%%%%%%%%%%%%%%%%%%%%%%%%%%%%%%%%%%%%%%%%%%%%%%%%%%%%%%%
Sandi Klav\v{z}ar was supported by the Slovenian Research Agency (ARIS) under the grants P1-0297, N1-0355, and N1-0285.

%%%%%%%%%%%%%%%%%%%%%%%%%%%%
\section*{Declaration of interests}
%%%%%%%%%%%%%%%%%%%%%%%%%%%%

The authors declare that they have no known competing financial interests or personal relationships that could have appeared to influence the work reported in this paper.

%%%%%%%%%%%%%%%%%%%%%%%%%%%%
\section*{Data availability}
%%%%%%%%%%%%%%%%%%%%%%%%%%%%

Our manuscript has no associated data.

\end{document}